\documentclass{ifacconf}

\usepackage{graphicx}      
\usepackage{natbib}        
\usepackage{additionalSetting}
\begin{document}
\begin{frontmatter}

\title{Fast identification and stabilization of unknown linear systems} 


\author[First]{Dennis Gramlich} 
\author[First]{Christian Ebenbauer}

\address[First]{Chair of Intelligent Control Systems,
	RWTH-Aachen,
	D-52074 Aachen, Germany
	{\tt\small \{dennis.gramlich,christian.ebenbauer\}@ic.rwth-aachen.de}}

\begin{abstract}
	
	In the present work, a simple algorithm for stabilizing an unknown linear time-invariant system is proposed, assuming only that this system is stabilizable. The suggested algorithm is based on first performing a partial identification of the system and then stabilizing the controllable subsystem. It should be emphasized that our approach does not depend on any prior model knowledge and requires only a minimal number of samples.
	
\end{abstract}


\end{frontmatter}

\section{INTRODUCTION}
\label{sec:1}

The ability to generate a stabilizing controller for an unknown system as quickly as possible (online) is a key capability for many control applications.
For example, many reinforcement learning/adaptive control methods require an initial stabilizing controller to be able to learn safely (\cite{fazel2018global,kakade2020information,tu2019gap}).
There exists the field of online identification, where a dynamical system is first identified at runtime and then controlled \cite{ward2010benefits,morelli1998flight}. Furthermore, in adaptive control learning triggers that initiate the re-identification of a system as soon as significant deviations occur between the previous system model and the current data have recently been introduced (see \cite{schluter2020event,solowjow2020event}).
For all these cases it is crucial that the determination of a (new) model and the stabilization of the system are done as fast as possible, since these are online applications and, hence, instability can cause damage to the system. The case of non-controllable, stabilizable systems, to which we dedicate much attention in this work, is also of practical relevance. For example, lifted systems or systems which incorporate disturbance models are stabilizable, but not controllable.

The obvious minimum requirement for performing this online (identification) stabilization task is the stabilizability of the considered system. Indeed, for the case of a controllable linear system with state dimension $n$ and input dimension $m$ the straightforward strategy to first apply a persistently exciting input signal of order $n+1$ \cite{willems2005note}, then identify the system matrices $(\hat{A},\hat{B})$ which are consistent with the measured data and finally design a stabilizing controller gain for $(\hat{A},\hat{B})$ generates a stabilizing controller for the system. This strategy for obtaining the stabilizing controller would be in line with standard online identification strategies. Using persistently exciting input signals to identify the system behaviour is still being actively researched in the field of data driven/adaptive control. For example, the algorithm derived in \cite{van2021beyond} from Willems' lemma resembles our method and achieves the identification of a controllable system in $n+m$ steps. However, despite the possibility of studying stabilizable systems, this line of research \cite{van2021beyond,van2020willems,verhoek2021fundamental,ist:berberich21g} mostly studies controllable systems and the identification of the full system behaviour. The possibility of generalizing Willems' Lemma to the stabilizable case is e.g. considered in \cite{yu2021controllability} and addressed in Section~\ref{sec:4}. 

The online stabilization of unknown (linear) systems is also actively researched in Reinforcement Learning, where research is conducted on regret bounds \cite{dean2018regret,mania2019certainty,simchowitz2020naive,abbasi2011regret,chen2021black} and finite time stabilization \cite{kazem2018finite,faradonbeh2019randomized,dean2020sample} for unknown linear dynamical systems. While older works on LQ Reinforcement Learning with bounded regret requires either an initial stabilizing controller or a controllable system, Reinforcement Learning with bounded regret for stabilizable systems without an initially stabilizing controller has been established in \cite{lale2022reinforcement}. However, unlike our work, \cite{lale2022reinforcement} relies on stochastic system noise to obtain sufficiently informative data.

As we see, the online stabilization and identification of unknown linear systems has been addressed many times, but 
none of the above mentioned works achieves stabilization and partial identification under minimal assumptions and in minimal time. The universal regulator \cite{maartensson1985order} as a prominent example for the stabilization of an unknown stabilizable system is neither fast nor simple. Therefore, we want to deal with exactly this question, i.e., to compute an asymptotically stabilizing feedback gain for an unknown stabilizable linear system as fast as possible, in the present article. In the course of our exposition, we explore the limits for the fastest and simplest, but not necessarily most robust stabilization of a dynamical system. We develop an elementary algorithm that achieves the partial identification and stabilization of a linear system in the minimum number of time steps. We emphasize that our algorithm does not rely on a persistently exciting input signal. Further, it does not perform a full system identification but rather a partial identification in the sense that it identifies that part of the system dynamics which is required for stabilization.

Contributions:
\begin{itemize}
	\item Our key contribution is a simple online algorithm (Algorithm \ref{alg:1}) for the fastest possible partial identification and stabilization of a stabilizable linear system, which does not rely on persistency of excitation.
	\item In addition to that, we provide a simple condition to be verified online (Theorem \ref{thm:4}) to ensure that any algorithm explores the full controllable subspace.
	\item Finally, we highlight that we work with the minimal assumptions for stabilization, which is stabilizability of the matrix pair $(A,B)$. 
\end{itemize}
\section{Problem statement}
\label{sec:2}

Consider a discrete-time linear dynamical system in state-space representation
\begin{align}
	x_{k+1} = A x_k + Bu_k, \label{eq:sys}
\end{align}
where the pair of matrices $(A,B) \in \bbR^{n\times n}\times \bbR^{n\times m}$ is completely unknown. We assume that the state $x_k \in \bbR^n$ can be measured, and denote the input by $u_k \in \bbR^m$.

The problem studied in this article is to find a stabilizing controller online. This, we formulate as follows.

\begin{problem}[Stabilization of an unknown system]~\\
	Given an unknown stabilizable pair of matrices $(A,B)$ and initial condition $x_0$. Choose (online) a sequence of inputs $(u_k)_{k=0}^{N-1} \subsEq \bbR^m, N \in \bbN$ such that a stabilizing controller $u_k = K x_k$ for \eqref{eq:sys} can be determined. This makes our setting more challenging than the case of controllable systems, where a complete identification of the system behaviour renders the stabilization obvious.
\end{problem}

\section{A simple solution}
\label{sec:3}

The simple solution we propose consists of first generating a data set of the form
\begin{align*}
	\calD_N := \{ (x_k^+,x_k,u_k)\mid k = 0,\ldots,N-1 \},
\end{align*}
where $x_k^+ = Ax_k + B u_k$, through interactions with the system (appropriate choice of $u_k$), then determining a pair of matrices $(\hat{A},\hat{B})$ which is stabilizable and consistent with $\calD_N$, and finally choosing a controller $K$, which stabilizes $(\hat{A},\hat{B})$ and therefore (Theorem~\ref{thm:2}) also $(A,B)$. 

One possible way to generate the required data would be to just apply a persistently exciting input signal $(u_k)_{k=0}^{N-1}$ of sufficient length $N$ to the system and to extract the required information. In this article, however, we consider a simpler solution which does not (directly) rely on persistency of excitation.
This simple solution is based on the following idea. If $x_0$ excites all modes of \eqref{eq:sys}, then we can just wait $n$ time steps (with $u_k = 0$) and then compute $\hat{A} = A$ by solving the equation
\begin{align*}
	\hat{A} = \begin{pmatrix}
		x_{n} & \cdots & x_1
	\end{pmatrix}
	\begin{pmatrix}
		x_{n-1} & \cdots & x_0
	\end{pmatrix}^{-1}
\end{align*}
if the matrix inverse exists. Thereafter, we can incrementally identify the columns of the matrix $B$ from \eqref{eq:sys} by setting $u_k$ equal to the $i$-th unit vector $e_i$  at time step $k = n+i$. Then we obtain
\begin{align*}
	b^{(i)} = x_{k+1} - A x_k,
\end{align*}
where we know all quantities on the right side of the equation and can therefore compute the $i$-th column $b^{(i)}$ of $B$. In this fashion, we can also identify $B$ after time step $k = n+m$. This extremely simple algorithm works for identifying $(A,B)$ in the case where $x_{n-1}, \cdots, x_0$ are linearly independent. This is satisfied with probability one, if $x_0$ is chosen randomly and independent of $A$ from a distribution which is absolutely continuous with respect to the Lebesque measure and if the kernel of $A$ has at most dimension $1$. Notice, that if it works, this algorithm identifies $(A,B)$ after $k = n+m$ steps, which is optimal.

However, in the present work, we do not want to rely on stochastic assumptions but instead, we want an algorithm capable of achieving the best possible identification of $(A,B)$ irrespective of the choice of the initial condition $x_0$ or the eigenvalues of $A$. 

The key is to observe that linear dependence of $x_0,\ldots,x_{n-1}$ can be detected as soon as $x_k \in \Span \{ 0,x_0, \ldots, x_{k-1}\}$ occurs, which we can easily check at any time step. Once $x_{k} \in \Span \{ 0,x_0, \ldots, x_{k-1}\}$ is detected, we know that the state will remain in this subspace forever without exterior excitation and we can also compute $A x_k$. Hence, in this case, we can directly switch $u_k = e_1$ and compute $b^{(1)} = x_{k+1} - Ax_k$ in the next time step and at the same time performing an exploration step in the direction of $b^{(1)}$. Now there are two possibilities: Either $x_{k+1} \notin \Span \{ 0,x_0, \ldots, x_{k}\}$ is observed, which means that we can go on identifying further information about $A$, or $x_{k+1} \in \Span \{ 0,x_0, \ldots, x_{k}\}$, which means that we can identify another column of $B$. This procedure is described in Algorithm \ref{alg:1}. At the end of this procedure, we have the data to identify the matrix $B$ completely and to identify $A$ on the smallest invariant subspace of $A$ which contains the vectors $\{x_0,b^{(1)},\ldots,b^{(m)}\}$. Achieving more is not possible with an online data set.

\begin{algorithm}
	\caption{Fast exploration}\label{alg:1}
	\begin{algorithmic}
		\STATE $x_{-1} \gets 0$,~~ $s \gets 1$,~~ $k \gets 0$,~~ $\calD_0 \gets \emptyset$
		\WHILE{true}
		\STATE $u_{k} \gets 0$
		\IF{$x_{k} \in \Span \{ x_{k-1},\ldots,x_{-1} \}$}
		\IF{$s> m$}
		\STATE \textbf{break}
		\ELSE
		\STATE $u_{k} \gets e_s$
		\STATE $s \gets s + 1$
		\ENDIF
		\ENDIF
		\STATE $x_{k+1} = A x_k + B u_k$ (interaction with the true sys.)
		\STATE $\calD_{k+1} = \calD_k \cup \{(x_{k+1},x_k,u_k)\}$
		\STATE $k \gets k+1$
		\ENDWHILE
		\RETURN $\calD_k$
	\end{algorithmic}
\end{algorithm}

Algorithm \ref{alg:1} provides us with a data set $\calD_N$, which means that we still need to do the (partial) identification of $(\hat{A},\hat{B})$ and to design a stabilizing controller gain $K$.

To this end, we define, based on a given data set $\calD_N = \{ (x_k,u_k,x_k^+)\mid k = 0,\ldots,N-1\}$, the data matrices
\begin{align*}
	X &:= \begin{pmatrix}
		x_0 & \cdots & x_{N-1}
	\end{pmatrix},\\
	X^+ &:= \begin{pmatrix}
		x_0^+ & \cdots & x_{N-1}^+
	\end{pmatrix}
	= \begin{pmatrix}
		x_1 & \cdots & x_{N-1} & x_{N-1}^+
	\end{pmatrix},\\
	U &:= \begin{pmatrix}
		u_0 & \cdots & u_{N-1}
	\end{pmatrix}.
\end{align*}
With the aid of these matrices, we can define the \emph{pseudo estimator} for $\hat{A}_0(\calD)$ and $\hat{B}_0(\calD)$ as follows.

\begin{definition}[Pseudo estimator]
	Consider a data set $\calD_N = \{ (x_k^+,x_k,u_k)_{k=0}^{N-1} \}$. We define the pseudo estimator $(\hat{A}_0(\calD_N),\hat{B}_0(\calD_N))$ as
	\begin{align*}
		(\hat{A}_0(\calD_N),	\hat{B}_0(\calD_N)) = \argmin_{\substack{(\hat{A},\hat{B}) \in \bbR^n \times \bbR^m,\\
		\forall (x_k^+,x_k,u_k) \in \calD_N :\\
		x^+_k = \hat{A} x_k + \hat{B} u_k}
				} \|\hat{A}\|_F^2 + \|\hat{B}\|_F^2,
	\end{align*}
	where $\|\cdot\|_F$ denotes the Frobenius norm or equivalently
	\begin{align*}
		\begin{pmatrix}
			\hat{A}_0(\calD_N) &
			\hat{B}_0(\calD_N)
		\end{pmatrix}
		=
		X^+
		\begin{pmatrix}
			X\\
			U
		\end{pmatrix}^\dagger,
	\end{align*}
	where $(\cdot)^\dagger$ denotes the Moore Penrose Pseudo inverse.
\end{definition}

For our purposes, the key property of the pseudo estimator $(\hat{A}_0(\calD_N),\hat{B}_0(\calD_N))$ is that it identifies a pair of matrices, which is \emph{consistent with the data} in $\calD_N$, i.e., $x^+ = \hat{A}_0(\calD_N)x + \hat{B}_0(\calD_N) u$ for any $(x^+,x,u) \in \calD_N$ and that it will set the eigenvalues of $\hat{A}$ belonging to unknown modes to zero, if we only have data from an invariant subspace of $A$ (which can happen). This property ensures that $(\hat{A}_0(\calD_N),\hat{B}_0(\calD_N))$ is going to be stabilizable, provided that $\calD_N$ is generated by Algorithm \ref{alg:1} and $(A,B)$ is stabilizable.

\begin{theorem}
	\label{thm:2} If $(A,B)$ is stabilizable and $\calD_N$ is a data set determined by Algorithm \ref{alg:1}, then also the pseudo estimate $(\hat{A}_0(\calD_{N}),\hat{B}_0(\calD_{N}))$ is stabilizable and any controller gain $K$ that stabilizes $(\hat{A}_0(\calD_{N}),\hat{B}_0(\calD_{N}))$ also stabilizes $(A,B)$.
\end{theorem}

In addition, the following theorem shows that we recover $(A,B)$ by $(\hat{A}_0(\calD_N),\hat{B}_0(\calD_N))$ if this is possible by online identification.

\begin{theorem}
	\label{thm:1} Algorithm \ref{alg:1} terminates after at most $n+m$ time steps. In addition, if there exists an $N \in \bbN$ and an input sequence $(u_k)_{k=0}^{N-1} \subsEq \bbR^m$ such that $(A,B)$ is uniquely determined (in the sense of Definition \ref{def:H}) by the resulting data set $\calD_N = \{ (x_{k+1},x_k,u_k) \mid k = 0,\ldots,N-1 \}$, then $(A,B)$ is also uniquely determined by $\calD_{N}$ from Algorithm~\ref{alg:1} and $(A,B) = (\hat{A}_0(\calD_{N}),\hat{B}_0(\calD_{N}))$.
\end{theorem}

The proofs of Theorem \ref{thm:1} and Theorem \ref{thm:2} follow from Theorem~\ref{thm:3} and can be found at the end of Section \ref{sec:3}.

\subsection{Properties of Algorithm \ref{alg:1}}

In the following, we view the identification task as the task to find a set of matrices (or a single matrix tuple), which is able to explain the data $\calD_N$. This generalization is necessary, because in the case of a pair of matrices $(A,B)$ which is not controllable, but only stabilizable, it might be impossible to uniquely determine $(A,B)$.

\begin{definition}[Set of data consistent parameters]
	\label{def:H}
	Let a dataset $\calD_N = \{(x_{k}^+,x_k,u_k) \mid k=0,\ldots, N-1\}$ be given. We call any pair of matrices $(\hat{A},\hat{B}) \in \bbR^{n\times n} \times \bbR^{n\times m}$ which is able to explain the data, i.e., which satisfies
	\begin{align*}
		x_{k}^+ &= \hat{A} x_k + \hat{B} u_k & k = 0,\ldots, N-1,
	\end{align*}
	consistent with the data $\calD_N$. The set of all matrices $(\hat{A},\hat{B})$ which are consistent with the data is the affine subspace $\calH(\calD_N)$ of $\bbR^{n\times n}\times \bbR^{n \times m}$. We say that $\calD_N$ uniquely determines $(A,B)$ if $\calH(\calD_N) = \{ (A,B) \}$.
\end{definition}

Ideally, the set of data consistent parameters $\calH(\calD_N)$ contains only one parameter $(\hat{A},\hat{B}) = (A,B)$. Unfortunately, this is not always the case. To ensure this, we would have to require that $(x_k,u_k)_{k=0}^{N-1}$ forms a basis of $\bbR^n \times \bbR^m$. This can easily be ensured when constructing an arbitrary data set (where the data must not be generated from a single trajectory), but it is a non-trivial problem when generating an online data set. Here, it is possible that not all of $\bbR^n \times \bbR^m$ can be explored. Indeed, depending on the initial condition $x_0$, we are able to define a subspace of $\bbR^n$, which can never be left by $(x_k)$.

\begin{definition}[Explorable subspace]
	We call a subspace $\calV \subsEq \bbR^n$ invariant with respect to $A$ if $A \calV \subsEq \calV$ and denote the set of invariant subspaces of $A$ with $\calI(A)$. The explorable subspace is the smallest invariant subspace of $A$ containing $\{x_0,b^{(1)},\ldots b^{(m)}\}$
	\begin{align*}
		\calV_{\text{exp}}(x_0) := \bigcap \{ \calV \in \calI(A) \mid \{x_0 , b^{(1)},\ldots, b^{(m)} \} \subsEq \calV \}.
	\end{align*}
\end{definition}

Clearly, any sequence $(x_k)$ defined by $x_{k+1} = A x_k + B u_k$ for some $(u_k)\subsEq \bbR^m, x_0 \in \bbR^n$ can never leave $\calV_{\text{exp}}(x_0)$ and, consequently, $(x_k^+,x_k,u_k) \in \calD_N$ for an online data set implies $x_{k}^+,x_k \in \calV_{\text{exp}}(x_0)$. Hence, with an online data set, we can at most hope to gather information about $A$ on $\calV_{\text{exp}}(x_0)$ (as preimage). In the light of these considerations, we can define the optimal set of data consistent matrices.

\begin{definition}[Optimal identifyable set]
	\label{def:4} Let some initial condition $x_0$ be given. We then define
	\begin{align*}
		\calH^*(x_0) := \{ (\hat{A},\hat{B}) \in \bbR^n\times \bbR^m \mid \hat{A} x + \hat{B} u = Ax + Bu,\\ (x,u)\in \calV_{\text{exp}}(x_0)\times \bbR^m\}.
	\end{align*}
\end{definition}

Any matrix from $\calH^*(x_0)$ is able to explain any online data set $\calD_N$ that could possibly be generated from $x_0$ (with any input sequence $(u_k)$). Hence, the best result that an exploration strategy could achieve in terms of identification is determining the set $\calH^*(x_0)$.

This identification of $\calH^*(x_0)$ is indeed possible. The following result shows that Algorithm \ref{alg:1} is able to identify $\calH^*(x_0)$ in minimum time.

\begin{theorem}
	\label{thm:3}
	Let $\calD_N$ be generated by Algorithm \ref{alg:1}. Then $\calH(\calD_N) = \calH^*(x_0)$ and $|\calD_N| = \tilde{n} + m \leq n + m$, where $\tilde{n} = \dim \calV_{\text{exp}}(x_0)$.
\end{theorem}

\begin{proof}
	Let $N$ be equal to the time index when Algorithm~\ref{alg:1} terminates, consider the data set $\calD_N$ returned by Algorithm~\ref{alg:1} and form the matrix
	\begin{align}
		\begin{pmatrix}
			X\\
			U
		\end{pmatrix}
		=
		\begin{pmatrix}
			x_0 & \cdots & x_{N-1}\\
			u_0 & \cdots & u_{N-1}
		\end{pmatrix} \label{eq:dataMatrix}
	\end{align}
	from its elements. Then the columns of this matrix are linearly independent, since every time when an entry $x_k$ depends linearly on the previous $x_0,\ldots, x_{k-1}$, the input $u_k$ is set to one of the unit vectors $e_1,\ldots, e_m$, which is linearly independent from the previous inputs $u_0,\ldots,u_{k-1}$. When linear dependence of $x_k$ occurs for the $m+1$-st time, the algorithm terminates. Note that this also implies that the matrix $X$ has rank $N-m$, since $m$ linearly dependent vectors have been added to $X$ during the runtime of Algorithm \ref{alg:1}. Next, we can permute the columns of \eqref{eq:dataMatrix} to obtain the matrix
	\begin{align*}
		\left(
		\begin{array}{ccc|ccc}
			x_{i_1} & \cdots & x_{i_{N-m}} & x_{i_{N-m+1}} & \cdots & x_{i_N}\\\hline
			0 & \cdots & 0 & e_1 & \cdots & e_m
		\end{array}
		\right).
	\end{align*}
	Since this matrix has full rank, both blocks on the diagonal must have full rank. Further, since $X$ has rank $N-m$, the vectors $x_{i_1},\ldots , x_{i_{N-m}}$ span the image of $X$. We can utilize this observation to show that $\range X$ is invariant w.r.t. $A$. To this end, it is sufficient to show $A x_{k} \in \range X$ for $k = i_1,\ldots,i_{N-m+1}$. There are two possibilities: If $k < N-1$, then $x_k^+ = Ax_k$, since $u_k = 0$ and $x_k^+ = x_{k+1}$ is one of the columns of $X$. If $k = N-1$, then Algorithm \ref{alg:1} has terminated at this index, which implies that $x_k^+ = x_N = Ax_k$ has been linearly dependent on the other columns $x_{k-1}, \ldots ,x_0$. Hence, $\range X$ is invariant w.r.t. $A$. It can also be shown that $\range X$ contains $\{x_0 ,b^{(1)},\ldots, b^{(m)}\}$. This is trivial for $x_0$. For $b^{(s)}$ and arbitrary $s$, we obtain
	\begin{align*}
		b^{(s)} = B e_s = A x_{i_{N-m+s}} - x_{i_{N-m+s}},
	\end{align*}
	where the right side is in $\range X$. Hence, since $\range X$ is obtained from a single trajectory, it must be $\range X = \calV_{\text{exp}}(x_0)$. Since already the first $N-m$ vectors form a basis of $\range X$, this implies that the matrix \eqref{eq:dataMatrix} spans $\calV_{\text{exp}}(x_0) \times \bbR^m$.
	
	The set $\calH^*(x_0)$ is defined as the set of solutions of the infinite system of equations
	\begin{align*}
		Ax + Bu &= \hat{A} x + \hat{B} u & x \in \calV_{\text{exp}}(x_0), u \in \bbR^m.
	\end{align*}
	in the decision variable $(\hat{A},\hat{B})$. The set $\calH(\calD_N)$ is defined as the solution set of the system of equations
	\begin{align*}
		x_k^+ &= A x_k + Bu_k = \hat{A} x_k + \hat{B} u_k & k = 0,\ldots,N-1
	\end{align*}
	in the variable $(\hat{A},\hat{B})$. Since the vectors $\begin{pmatrix}
		x_k^\top & u_k^\top
	\end{pmatrix}^\top$ span $\calV_{\text{exp}}(x_0) \times \bbR^m$, these solution sets are the same.
	
	Finally, since the columns of \eqref{eq:dataMatrix} are linearly independent and their span $\calV_{\text{exp}}(x_0) \times \bbR^m$ has dimension $\tilde{n} + m$, it must be $N = \tilde{n} + m$.
\end{proof}

\emph{Proof of Theorem \ref{thm:1}}: The termination after $k \leq n+m$ follows trivially from Theorem \ref{thm:3}, because $\tilde{n} + m \leq n + m$. In addition, $(A,B)$ can only be determined uniquely with some input sequence, if $\calH^*(x_0) = \{(A,B)\}$. However, Algorithm \ref{alg:1} constructs $\calH^*(x_0)$ and therefore determines $(A,B)$ whenever there is an input sequence which can achieve this identification. \hfill $\blacksquare$

To prove Theorem \ref{thm:2}, it is important to understand that two pairs of matrices with the same controllable subspace and the same dynamics on this controllable subspace are stabilized by the same controller gains.

\begin{lemma}
	\label{lem:controllableSubspace}
	Consider $(A^1,B^1), (A^2,B^2) \in \bbR^{n\times n}\times \bbR^{n\times m}$ and let $\range B^1 \cup \range B^2 \subsEq \calV \subsEq \bbR^n$ be an invariant subspace w.r.t. $A^1$ and $A^2$ such that $A^1 x + B^1 u = A^2 x + B^2 u$ for all $(x,u) \in \calV\times \bbR^m$. Then $A^1 + B^1 K$ is stable if and only if $A^2 + B^2 K$ is stable.
\end{lemma}
\begin{proof}
	Consider a state transform
	\begin{align}
		x = T \tilde{x} = \begin{pmatrix}
			T_1 & T_2
		\end{pmatrix} \tilde{x}, \qquad T^\top T = I, \label{eq:tildeCoordinates}
	\end{align}
	where $T_1$ forms a basis matrix of $\calV$. In the new coordinates $\tilde{x}$, $A^{1} + B^1 K$ and $A^2 + B^2 K$ take the form
	\begin{align*}
		\begin{pmatrix}
			\tilde{x}^{(1)}_{k+1}\\
			\tilde{x}^{(2)}_{k+1}
		\end{pmatrix}
		=
		\left(
		\begin{pmatrix}
			A_{11}^{i} & A_{12}^{i}\\
			0 & A_{22}^{i}
		\end{pmatrix}
		+
		\begin{pmatrix}
			B_1^i\\
			0
		\end{pmatrix}
		\begin{pmatrix}
			K_1 & K_2
		\end{pmatrix}
		\right)
		\begin{pmatrix}
			\tilde{x}^{(1)}_{k}\\
			\tilde{x}^{(2)}_{k}
		\end{pmatrix},
	\end{align*}
	because $\calV$ is invariant w.r.t. $A^{i}$, $i = 1,2$, and the columns of $B^i$ are contained in $\calV$. The stability of these systems only depends on the stability of $A_{11}^{i} + B_1^i K_1$, since  $A_{22}^{1}$ and $A_{22}^{2}$ must be stable due to the stabilizability assumption on $(A^{1},B^1)$ and $(A^{2},B^2)$. However, $A_{11}^{1} + B_1^1 K_1$ and $A_{11}^{2} + B_1^2 K_1$ are equal, because $(A^{1},B^2)$ and $(A^{2},B^2)$ are equal on $\calV\times \bbR^m$. Consequently, $K$ stabilizes $(A^{1},B^1)$ if and only if it stabilizes $(A^{2},B^2)$.
\end{proof}

Due to Lemma \ref{lem:controllableSubspace}, a stabilization scheme must only make sure to determine $(A,B)$ on the controllable subspace, i.e., the smallest invariant $\calV$ such that $\range B \subsEq \calV$. This conclusion is the key to prove Theorem \ref{thm:2}.

\emph{Proof of Theorem \ref{thm:2}}: According to Theorem \ref{thm:3}, $\calH (\calD_N) = \calH^*(x_0)$ holds for the data set $\calD_N$ returned by Algorithm~\ref{alg:1}. Consequently, all matrix pairs in $\calH(\calD_N)$ are equal on $\calV_{\text{exp}}(x_0) \times \bbR^m$. Hence, by Lemma \ref{lem:controllableSubspace}, a controller which stabilizes some element of $\calH(\calD_N)$ stabilizes all stabilizable matrix pairs in $\calH(\calD_N)$ and all that is left to show is that $(\hat{A}(\calD_N),\hat{B}(\calD_N))$ extracts a stabilizable matrix pair from $\calH(\calD_N)$. To see that the latter is true, consider that $\calH(\calD_N)$ takes the form
\begin{align*}
	\resizebox{\linewidth}{!}{$
		\left\{ \left( \begin{pmatrix}
			A_{11} & \hat{A}_{12}\\
			0 & \hat{A}_{22}
		\end{pmatrix}, \begin{pmatrix}
			B_1\\
			0
		\end{pmatrix} \right) \mid \hat{A}_{22} \in \bbR^{(n-\tilde{n})\times (n-\tilde{n})}, \hat{A}_{12} \in \bbR^{\tilde{n} \times (n-\tilde{n})} \right\},$}
\end{align*}
in the coordinates \eqref{eq:tildeCoordinates}, where $\hat{A}_{11} = A_{11}$ and $\hat{B}_1 = B_1$ as these matrices are determined by the data and $\hat{A}_{12}$ and $\hat{A}_{22}$ are completely free. Since the pseudo estimator minimizes the Frobenius norm among all elements of $\calH(\calD_N)$, it will estimate $\hat{A}_{12}$ and $\hat{A}_{22}$ to zero and thus extract a stabilizable matrix pair. Since the Frobenius Norm, and therefore also the pseudo estimator, is invariant w.r.t. orthogonal coordinate changes, this result generalizes to all coordinates.
\hfill $\blacksquare$

\section{A charaterization for fast exploration}
\label{sec:4}

A major challenge of indirect adaptive control, i.e., control methods that split the stabilization of the system into a system identification and a controller synthesis, is to extract stabilizable models from the data. Once this problem is solved, a controller can be generated which either stabilizes the system or explores a deviation between the previous identification and the true system after at most $n$ steps. Unfortunately, the problem of finding a stabilizable model that is consistent with a data set $\calD_N$ is generally non-convex. However, in a perturbation-free setting, convex extraction of a data consistent stabilizable system can be accomplished by appropriate exploration. Algorithm \ref{alg:1} shows that this is possible in $\tilde{n} + m$ steps, and explores a perturbation-free system in this minimal number of steps in the best possible way. The next theorem presents a condition which is sufficient for these desirable properties of Algorithm \ref{alg:1}. In the further course of this section, we also compare our exploration strategy to persistently exciting input signals as used in the so-called Willems' Lemma \cite{willems2005note}. 

\begin{theorem}
	\label{thm:4}
	Assume we are given an online data set $\calD_N = \{ (x_k^+,x_k,u_k) \mid k = 0,\ldots , N-1 \}$ with $x_k^+ = Ax_k + Bu_k~ \forall k = 0,\ldots,N-1$ and $x_k = x_{k-1}^+$. Then the following two statements hold true.
	\begin{enumerate}
		\item $\begin{pmatrix}
			u_{N}\\
			x_{N-1}^+
		\end{pmatrix}
		\notin \Span \left\{ \begin{pmatrix}
			u_0\\
			x_0
		\end{pmatrix}, \ldots ,\begin{pmatrix}
			u_{N-1}\\
			x_{N-1}
		\end{pmatrix} \right\}$ implies strict inclusion $\calH(\calD_{N+1}) \subset \calH(\calD_N)$.
		\item $\begin{pmatrix}
			\tilde{u}\\
			x_{N-1}^+
		\end{pmatrix}
		\in \Span \left\{ \begin{pmatrix}
			u_0\\
			x_0
		\end{pmatrix}, \ldots ,\begin{pmatrix}
			u_{N-1}\\
			x_{N-1}
		\end{pmatrix} \right\} \quad \forall \tilde{u} \in \bbR^m$ if and only if $\calH (\calD_N) = \calH^* (x_0)$.
	\end{enumerate}
\end{theorem}
\begin{proof}
	We start by showing 1). The condition
	\begin{align*}
		\begin{pmatrix}
			u_{N}\\
			x_{N-1}^+
		\end{pmatrix}
		\notin \Span \left\{ \begin{pmatrix}
			u_0\\
			x_0
		\end{pmatrix}, \ldots ,\begin{pmatrix}
			u_{N-1}\\
			x_{N-1}
		\end{pmatrix} \right\}
	\end{align*}
	implies that the set
	\begin{align*}
		\calD_{N+1} 
		= \calD_N \cup 
		\{ (Ax_N^+ + B u_{N+1},x_N^+,u_{N+1}) \}
	\end{align*}
	adds an additional linearly independent constraint on $(\hat{A},\hat{B})$ to $\calD_N$. Consequently, $\calH(\calD_{N+1}) \subset \calH(\calD_N)$.
	
	Now, we can show statement 2). Hence, we first assume
	\begin{align}
		\begin{pmatrix}
			\tilde{u}\\
			x_{N-1}^+
		\end{pmatrix}
		\in \Span \left\{ \begin{pmatrix}
			u_{N-1}\\
			x_{N-1}
		\end{pmatrix},
		\ldots ,\begin{pmatrix}
			u_0\\
			x_0
		\end{pmatrix}
		\right\} ~ \forall \tilde{u} \in \bbR^m \label{eq:inv}
	\end{align}
	is satisfied at time $N$. We show that this implies
	\begin{align*}
		\begin{pmatrix}
			u\\
			x
		\end{pmatrix}
		\in \calS_N := \Span \left\{ \begin{pmatrix}
			u_{N-1}\\
			x_{N-1}
		\end{pmatrix},
		\ldots ,\begin{pmatrix}
			u_0\\
			x_0
		\end{pmatrix}
		\right\}
	\end{align*}
	for any $x \in \Span \{x_0,\ldots,x_{N-1}\}$ and $u \in \bbR^m$. To this end, use \eqref{eq:inv} to choose $\lambda_0,\ldots \lambda_{N-1}$, such that $x_{N-1}^+ - x = \lambda_0 x_0 + \ldots + \lambda_{N-1} x_{N-1}$ and choose $u_{N} \in \bbR^m$ such that
	\begin{align*}
		\underbrace{
			\begin{pmatrix}
				u_{N}\\
				x_{N-1}^+
		\end{pmatrix}}_{\in \calS_N}
		=
		\begin{pmatrix}
			u\\
			x
		\end{pmatrix}
		+
		\underbrace{
			\lambda_{N-1}
			\begin{pmatrix}
				u_{N-1}\\
				x_{N-1}
			\end{pmatrix}
			+
			\ldots
			+
			\lambda_0
			\begin{pmatrix}
				u_0\\
				x_0
		\end{pmatrix}}_{\in \calS_N}.
	\end{align*}
	Hence, any state input pair we could choose with an initial state in $\calD_{N}$ is already contained in $\calS_N$. Further, due to \eqref{eq:inv}, $\calS_N$ is invariant, i.e., it cannot be left by $(x,u)$. Consequently, $\calS_N$ also containes any dataset generated by Algorithm \ref{alg:1}, which implies $\calH(\calD_N) = \calH^*(x_0)$. Finally, notice that if \eqref{eq:inv} were not satisfied, then 1) implies that $\calH(\calD_N)$ can be improved, which contradicts $\calH(\calD_N) = \calH^*(x_0)$. Hence, the converse direction holds as well.
\end{proof}
Theorem \ref{thm:4} describes the essential property of Algorithm~\ref{alg:1} that enables fast exploration. Specifically, the inputs chosen by Algorithm \ref{alg:1} always satisfy condition 1), so that $\calH(\calD_k)$ is more tightly constrained at each step $k$. Once 2) is satisfied, Algorithm \ref{alg:1} terminates. Condition~1) from Theorem \ref{thm:4} is also proposed in \cite{van2021beyond} for exploration methods of controllable linear systems. It follows from Theorem \ref{thm:4} that any exploration scheme satisfying 1) inherits the result of Theorem \ref{thm:3} from Algorithm \ref{alg:1}.

An alternative to the exploration strategy we propose is to use a persistently exciting input signal $(u_k)_{k=0}^{N-1}$ of order $L = n+1$ according to the following definition.

\begin{definition}[Persistency of excitation]
	An input signal $(u_k)_{k=0}^{N-1}$ is said to be persistently exciting of order $L$, if the Hankel matrix
	\begin{align*}
		H_L((u_k)_{k=0}^{N-1})
		:=
		\begin{pmatrix}
			u_0 & u_1 & \cdots & u_{N-L}\\
			\vdots & \vdots & \ddots & \vdots\\
			u_{L-1} & u_L & \cdots & u_{N-1}
		\end{pmatrix}
	\end{align*}
	has full row rank.
\end{definition}

This notion of persistency of excitation is common in the literature on Willems' Lemma. The following corollary, which follows from the generalized Willems' Lemma (Theorem 1) in \cite{yu2021controllability}, states that the order $L = n+1$ is sufficient for $\calH(\calD_N) = \calH^*(x_0)$. Note that there is a major difference to the result achieved by Algorithm \ref{alg:1}, since we consider now the case, where the input signal $(u_k)$ is fixed from time $k = 0$, whereas Algorithm \ref{alg:1} chooses the input signal online and makes use of information which is only online available.

\begin{corollary}
	\label{cor:aPrioriExploration}
	If a given input signal $(u_k)_{k=0}^{N-1} \subsEq \bbR^m$ is persistently exciting of order $L = n+1$, then it achieves $\calH(\calD_{N}) = \calH^*(x_0)$ for any matrices $(A,B) \in \bbR^{n\times n} \times \bbR^{n\times m}$.
\end{corollary}

We mention this result, firstly, to compare our approach with existing work in the literature and, secondly, because Corollary \ref{cor:aPrioriExploration} provides an interesting alternative to Algorithm \ref{alg:1} with advantages and disadvantages.

For instance, it is a clear advantage of Corollary \ref{cor:aPrioriExploration} that the input signal $(u_k)$ can be chosen a priori, i.e., it must not be chosen online, since persistency of excitation depends in no way on $(A,B)$ or the values $(x_0,\ldots,x_{N-1})$. On the other hand, however, one should note that a signal $(u_k)$, which is persistently exciting of order $(n+1)$ will be significantly longer than $\tilde{n}+m$. To see this, consider that the matrix $H_{n+1}((u_k)_{k=0}^{N-1})$ has $(n+1)m$ rows. Since it must have at least as many columns to have full row rank, we must require $N \geq 2(n+1)m$. This additional length is the price which is paid for the convenience of choosing the signal $(u_k)$ offline. Algorithm~\ref{alg:1} on the other hand can achieve faster identification, because it makes use of online knowledge for choosing $(u_k)$.

\section{Simulation Example}
\label{sec:5}

To check the viability of our proposed method, we conduct a simulation study using OpenAI Gym's Lunar Lander. We believe that the Lunar Lander simulates an interesting application, where a stabilizing controller needs to be determined quickly before the landing. For our method, this is a very challenging example, since the simulation is highly nonlinear and the task is to perform a landing and not just stabilization. The Lunar Lander has the six continuous states $(x,y,v_x,v_y,\psi,\omega)$, where $(x,y)$ is the position in two-dimensional space, $(v_x,v_y)$ is the velocity in two-dimensional space and $(\psi,\omega)$ is the orientation and angular velocity of the spacecraft. There are also two inputs $(u_1,u_2)$ which control the thrusters of the space craft. For detailed information on the Lunar Lander, we refer to \cite{brockman2016openai}. In order to identify the system, we perform the exploration with Algorithm \ref{alg:1} several times in some cases, since the system has perturbations due to nonlinearities and random initialization. The stochastic input perturbations in the system are particularly serious. Therefore, we consider identifications with input perturbations switched on and switched off. To simulate a landing with our stabilization technique, we specify a constant velocity reference trajectory for landing and run 1000 simulations for each of our scenarios. The simulation results can be found in Table \ref{tab:lunarLander}. Here we provide the percentage of successful landings and the average reward. In comparison with other learning algorithms used on the Lunar Lander, it should be noted that most procedures are trained over perhaps hundreds of episodes, while we consider the performance in the first episode.

\begin{table}
	\centering
	\caption{Rate of successful landings and average reward on the first episode when applying Algorithm \ref{alg:1} for the identification of the Lunar Lander. Here, the number $\ell$ refers to the number of calls of Algorithm \ref{alg:1} for the identification of the model.}
	\label{tab:lunarLander}
	\begin{tabular}{@{}llllll@{}}
		\toprule
		& Noise & $\ell = 1$ & $\ell = 2$ & $\ell = 3$
		\\\midrule
		success rate & No & 0.42 & 0.94 & 0.73\\
		avg. reward & No & -163 & 217 & 112\\
		success rate & Yes & 0.06 & 0.9 & 0.18\\
		avg. reward & Yes & -483 & -288 & -191
		\\\bottomrule
	\end{tabular}
\end{table}

The results in Table \ref{tab:lunarLander}, show that Algorithm \ref{alg:1} enables the emergency stabilization of dynamical systems. The low success rates in the noisy case indicate a lack of robustness of Algorithm \ref{alg:1}, which had to be expected in the limit of minimum data points. An interesting observation is that in the noise-free setting for $\ell = 3$ calls of Algorithm \ref{alg:1}, the success rate for landing is lower than in the $\ell = 2$ case, even though more data is used when $\ell = 3$. This phenomenon can probably be explained by the fact that in the case $\ell = 3$ the exploration takes longer, so that there might not be enough time for the emergency landing.


Our simulation is available on github under the link (\href{https://github.com/SphinxDG/FastStabilization}{https://github.com/SphinxDG/FastStabilization}). Note that this simulation also contains experiments with non-controllable linear systems, where it is confirmed that linear systems are stabilized by the proposed algorithm except for numerically caused exceptions.

\section{Conclusion}
\label{sec:6}

In the present work, we study the problem of stabilizing an arbitrary unknown stabilizable linear system. In terms of Algorithm \ref{alg:1}, we provide a simple solution, which solves the task after a minimum number of time steps. Algorithm~\ref{alg:1} is quite simple and can be extended in many ways. For example, the inputs $u_k$ must not be set to unit vectors in Algorithm \ref{alg:1}, when they are excited, but can be chosen as arbitrary linearly independent vectors, or an initial gain $K_0$ could be chosen based on prior knowledge, such that Algorithm \ref{alg:1} identifies $(A+BK_0,B)$ instead of $(A,B)$. Further research directions could investigate the robustness of Algorithm \ref{alg:1} or its integration into receding horizon control of unknown systems \cite{ebenbauer2021control}, but these are beyond the scope of this paper. Finally, we provide in Theorem \ref{thm:4} a simple condition to achieve the best possible exploration with any algorithm. Based on this result, we can e.g. pick the recursive least squares algorithm and perturb the input a little bit, when the condition in Theorem \ref{thm:4} is not met, to ensure that the best possible partial identification is achieved. 

\bibliographystyle{abbrv}
\bibliography{sources.bib}                                          
\end{document}